\documentclass[journal]{IEEEtran}
\DeclareUnicodeCharacter{2009}{\,}

\usepackage{orcidlink}
\usepackage{cite}

\usepackage{graphicx}
\usepackage{amsmath}
\usepackage{amsfonts}
\usepackage{amssymb}

\newcommand{\norm}[1]{\big\|#1\big\|}

\newcommand{\VEC}[1]{\mathbf{#1}}
\newcommand{\CORE}[2]{\mathcal{#1}_{#2}}

\newcommand{\TENS}[1]{\mathbf{#1}}
\newcommand{\TENSOP}[1]{\Hat{\mathbf{#1}}}

\newcommand{\En}[2]{\TENS{E}^{#1}_{#2}}

\newcommand{\Hn}[2]{\TENS{H}^{#1}_{#2}}
\newcommand{\HDotn}[2]{\Dot{\TENS{H}}^{#1}_{#2}}
\usepackage{tikz}

\newlength{\defaulttensorsize}
\newcommand{\tensorsize}{\defaulttensorsize}

\usetikzlibrary{shapes.geometric}
\usetikzlibrary{shapes.misc}
\usetikzlibrary{arrows.meta}

\tikzset{>=latex}

\tikzset{every picture/.style={line width=0.6pt}}

\tikzstyle{core}=[circle,fill=Accent1, minimum size=\tensorsize]
\tikzstyle{opcore} = [rectangle,fill=Accent1, minimum size=\tensorsize]

\tikzstyle{nocore}=[circle, minimum size=\tensorsize]

\tikzstyle{lorthcore}=[rounded rectangle, fill=Accent1, rounded rectangle left arc=none, minimum size=\tensorsize*1.3]
\tikzstyle{rorthcore}=[rounded rectangle, fill=Accent1, rounded rectangle right arc=none, minimum size=\tensorsize*1.3]
\tikzstyle{diagcore}=[diamond, fill=Accent2, , minimum size=0.8\tensorsize]

\tikzstyle{widetensor}[2]=[rectangle, rounded corners = 8pt, fill=Accent1, minimum width=\defaulttensorsize*#1,minimum height = \defaulttensorsize]

\tikzstyle{tensornetwork}=[baseline=-0.25em]

\usepackage{xcolor}

\definecolor{Accent2}{RGB}{235, 85, 77}
\definecolor{Main}{RGB}{18, 92, 93}

\colorlet{Accent1}{Main}

\usepackage{hyperref}
\usepackage[capitalize]{cleveref}
\crefname{equation}{}{}

\usepackage{standalone}
\usepackage{booktabs}

\usepackage{siunitx}

\newcommand\submittedtext{%
  \footnotesize\itshape This work has been submitted to the IEEE for possible publication.}

\newcommand\submittednotice{%
\begin{tikzpicture}[remember picture,overlay]
\node[anchor=south,yshift=10pt] at (current page.south) {\parbox{\dimexpr0.65\textwidth-\fboxsep-\fboxrule\relax}{\centering\submittedtext}};
\end{tikzpicture}%
}

\begin{document}
\markboth{Submitted to IEEE Antennas and Wireless Propagation Letters for possible publication.}{}
\title{A Highly Scalable Quantized Tensor-Train FDTD Framework for the Simulation of Three-Dimensional Electromagnetic Scattering Problems}
\author{Daan Vanhaecke\orcidlink{0009-0009-2727-9252}, Emile Vanderstraeten\orcidlink{0000-0001-5189-388X} \IEEEmembership{Member, IEEE}, Dries Vande Ginste\orcidlink{0000-0002-0178-288X} \IEEEmembership{Senior Member, IEEE}
\thanks{Received -- -------- ----; revised -- --------- ----; accepted -- ------ ----.
Date of publication -- -------- ----; date of current version 15 Septemner 2026. (Corresponding author: Daan Vanhaecke)}
\thanks{The authors are with quest, Department of Information Technology, Ghent University- imec, Belgium (e-mail: daavhaec.vanhaecke@ugent.be)}\thanks{Digital Object Identifier  xx}
}

\maketitle
\submittednotice
\begin{abstract}
    In this letter, a novel 3-D Finite-Difference Time-Domain (FDTD) framework is proposed that circumvents the costly volumetric scaling of conventional FDTD methods. By representing the electromagnetic fields as low-rank Quantized Tensor Trains (QTT), the memory requirements scale logarithmically with system size. Moreover, the various numerical operations that constitute the FDTD scheme can be efficiently implemented in this format, with their computational cost also exhibiting a logarithmic complexity. The simulation of systems in open space is enabled by the inclusion of a uniaxial PML.  A validation example demonstrates that the method achieves excellent accuracy compared to the traditional full-grid (FG) FDTD method, while significantly reducing the required computational resources. Memory savings of several orders of magnitude are obtained, highlighting the potential of the proposed framework for the simulation of large multiscale electromagnetic systems.
\end{abstract}

\section{Introduction}
\IEEEPARstart{A}{s} one of the most popular full-wave solvers the Finite-Difference Time-Domain (FDTD) technique \cite{yee_numerical_1966,taflove_computational_2005} is routinely employed in many fields, including microwave and optical devices, geophysics, plasma physics and medical imaging \cite{teixeira_finite-difference_2023}. Unfortunately, the method requires a volumetric discretization of the domain resulting in an $\mathcal{O}(n^3)$ memory requirement and cost per time step, where $n$ denotes the number of gridpoints along each dimension. Consequently, the simulation of electrically large systems quickly demands substantial resources. Even though the computational cost can be alleviated by parallelized and GPU execution, the problematic fundamental scaling behavior remains.

Recently, the use of low-rank tensor decompositions was proposed to reduce the complexity of the FDTD method. More specifically, instead of the conventional full-grid (FG) 3-D tensor format, the tensor-train (TT) decomposition was employed to represent the discretized electromagnetic fields, resulting in an $\mathcal{O}(n)$ cost per time step \cite{manzini_tensor-train_2023,zhou_tensor-train_2025,scherzer_rank-limiting_2026}. This letter, however, focuses on the \emph{quantized} version of the TT-decomposition, namely the quantized tensor train (QTT) decomposition. The QTT-decomposition additionally decomposes the data along a certain dimension using a quantized dyadic grid, in this way achieving additional compression. QTT-decompositions have been incorporated into a variety of Maxwell equation solvers, including a method of moments solver of volume integral equations \cite{nguyen_tensor_2026-1},  2-D Yee-FDTD solvers utilizing different types of time integration methods\cite{ye_practical_2026}, a 2-D $\text{TM}_{z}$-formulation of scattered field FDTD with Mur absorbing boundary conditions \cite{nguyen_tensor_2026}, 3-D vector potential based FDTD with support for several time-integrators \cite{ye_practical_2026} and a 3-D FDTD solver employing the related QTT-Tucker decomposition and Crank-Nickolson like time stepping \cite{ye_quantized_2024}. 

Whereas QTT-acceleration of explicit 3-D Yee-FDTD has been suggested \cite{nguyen_tensor_2026,scherzer_rank-limiting_2026}. This technique still remains to be fully investigated, developed and implemented. Additionally, more effective absorbing boundary conditions are desirable for the simulation of open-space scattering problems. Therefore, this letter presents the first fully three-dimensional QTT-accelerated explicit Yee-FDTD framework for electromagnetic scattering simulations in open space, including implementation and validation. The framework supports dielectric materials, utilizes step-and-truncate (SAT) time integration and employs an uniaxial perfectly matched layer (UPML) for domain truncation. The PML is implemented using auxiliary differential equations (ADEs).

The letter is structured as follows: \cref{sec:formulation} briefly introduces the basic notions of tensor trains and quantized tensor trains. Next, it provides a high-level overview of the framework and discusses the theoretical scaling of the novel scheme. In \cref{sec:numerical results} an illustrative validation example demonstrates the excellent accuracy and logarithmic scaling of the method. Finally, some concluding remarks are given in \cref{sec:conclusion}.

\section{Formulation}\label{sec:formulation}
\subsection{Tensor Trains}
The TT-decomposition of a $d$-dimensional tensor $\TENS{A}$ is defined by the element-wise equality
\begin{equation}\label{eq: tensor train}
        \TENS{A}(i_1,i_2,\ldots,i_d) = \CORE{A}{1}(i_1)\CORE{A}{2}(i_2)\ldots\CORE{A}{d}(i_d),
\end{equation}
where $\CORE{A}{k}(i_k)$ represents a $r_{k-1}\times r_{k}$ matrix $\forall i_k \in 1,\ldots, n_k$. The matrices $\CORE{A}{k}(i_k)$ are called the cores of the tensor train, $r_k$ the tensor train ranks and $n_k$ the mode sizes. Additionally, it holds that $r_0=r_d=1$. A TT is conveniently visualized using tensor network notation (TNN) as
\begin{equation}
    \scalebox{0.5}{
\begin{tikzpicture}[tensornetwork]
        \node[widetensor = 4] (A) at (0,0) {};

        \coordinate (A1) at (-1.5*\tensorsize, 0) {};
        \coordinate (A2) at (-0.5*\tensorsize, 0) {};
        \coordinate (A3) at (0.5*\tensorsize, 0) {};
        \coordinate (An) at (1.5*\tensorsize, 0) {};
        
        \draw (A.south -| A1) -- +(0, -0.5) node[below] {};
        \draw (A.south -| A2) -- +(0, -0.5) node[below] {};
        \draw (A.south -| A3)  +(0, -0.25) node[below] {$\cdots$};
        \draw (A.south -| An) -- +(0, -0.5) node[below] {};
       
\end{tikzpicture} 
} = \scalebox{0.5}{
\begin{tikzpicture}[tensornetwork]
    \node[core] (A1) at (0,0) {};
    \node[core] (A2) at (1.5,0) {};
    \node[nocore] (A3) at (3,0) {$\ldots$};
    \node[core] (An) at (4.5,0) {};

    \draw (A1) -- (A2) -- (A3) -- (An);
    \draw (A1.south) -- + (0,-0.5);
    \draw (A2.south) -- + (0,-0.5);
    \draw (An.south) -- + (0,-0.5);
\end{tikzpicture}
} \,.
\end{equation}
In TNN, tensors are depicted by geometrical shapes, their indices by lines emerging from them and contractions are represented by connected lines, see \cite{bridgeman_hand-waving_2017} for more details. The storage cost of a tensor in TT-format scales as $\mathcal{O}(dnr^2)$, where, for simplicity $r_k = r$ and $n_k = n$, $\forall k$, was assumed. When the rank $r$ remains moderate, this results in significant compression compared to the original storage and the problematic $\mathcal{O}(n^d)$ complexity is reduced to a linear scaling with respect to $n$. Moreover, in many cases the rank can be significantly reduced by resorting to an approximate TT-decomposition with some error tolerance $\varepsilon$ or fixed rank, which can be calculated using algorithms such as TT-SVD \cite{oseledets_tensor-train_2011} or TT-Cross \cite{oseledets_tt-cross_2010}. TT-SVD constructs the TT-representation by consecutive (truncated) singular value decompositions (SVD), while TT-Cross is a tensor sketching algorithm.
Linear operators acting on tensors can also be brought in the TT-format. They are represented by tensor train operators (TTO), which are very similar to tensor trains, except that their cores are indexed by two indices instead of one \cite{oseledets_approximation_2010_good,kazeev_low-rank_2012}. The TTO-decomposition of a linear operator $\mathcal{O}$ is given by the element-wise equality
\begin{equation}
    \TENSOP{O}(i_1,\ldots,i_d,j_1,\ldots,j_d) = \CORE{O}{1}(i_1,j_1)\ldots\CORE{O}{d}(i_d,j_d).
\end{equation}
Herein, the cores and ranks are defined completely analogously to the regular TT case. In TNN, TTOs are represented as
\begin{equation}
    \scalebox{0.5}{
\begin{tikzpicture}[tensornetwork]
        \node[widetensor = 4] (A) at (0,0) {};

        \coordinate (A1) at (-1.5*\tensorsize, 0) {};
        \coordinate (A2) at (-0.5*\tensorsize, 0) {};
        \coordinate (A3) at (0.5*\tensorsize, 0) {};
        \coordinate (An) at (1.5*\tensorsize, 0) {};
        
        \draw (A.south -| A1) -- +(0, -0.5) node[below] {};
        \draw (A.south -| A2) -- +(0, -0.5) node[below] {};
        \draw (A.south -| A3)  +(0, -0.25) node[below] {$\cdots$};
        \draw (A.south -| An) -- +(0, -0.5) node[below] {};

        \draw (A.north -| A1) -- +(0, 0.5) node[above] {};
        \draw (A.north -| A2) -- +(0, 0.5) node[above] {};
        \draw (A.north -| A3)  +(0, 0.25) node[above] {$\cdots$};
        \draw (A.north -| An) -- +(0, 0.5) node[above] {};
       
\end{tikzpicture} 
} = \scalebox{0.5}{
\begin{tikzpicture}[tensornetwork]
    \node[opcore] (O1) at (0,0) {};
    \node[opcore] (O2) at (1.5,0) {};
    \node[nocore] (O3) at (3,0) {$\ldots$};
    \node[opcore] (On) at (4.5,0) {};

    \draw (O1) -- (O2) -- (O3) -- (On);
    \draw (O1.south) -- + (0,-0.5);
    \draw (O2.south) -- + (0,-0.5);
    \draw (On.south) -- + (0,-0.5);

    \draw (O1.north) -- + (0,0.5);
    \draw (O2.north) -- + (0,0.5);
    \draw (On.north) -- + (0,0.5);

\end{tikzpicture}
}.
\end{equation}
\subsection{Quantized Tensor Trains}\label{sec: QTT}
Owing to the $\mathcal{O}(dnr^2)$ storage cost, higher compression rates can be attained by folding the tensor in such a way as to increase the dimensionality $d$ while decreasing $n$.  Thus, alternative to directly employing TTs, the dimensionality of the tensor can first be artificially increased by downfolding, after which a TT-decomposition is applied. This decomposition is referred to as the quantized tensor-train (QTT) decomposition \cite{khoromskij_odlog_2011,oseledets_approximation_2009}.

For example, a vector $\VEC{x}$ (order-one tensor) of length $n=2^l$ can be reshaped into an $l$-dimensional tensor $\TENS{X}$ with all mode sizes equal to $2$. In this work, a sequential ordering is adopted, meaning that the reshaped tensor $\TENS{X}$ is addressed by a multi-index $(i_1,\ldots,i_l)$, such that the original index $i$ is related to the new one via its binary string $i = (i_1,\ldots,i_l)_2$. The TT-decomposition is subsequently applied. For a general tensor, this downfolding strategy is adopted along each dimension, resulting in the decomposition depicted by
\begin{equation}
    \hbox{\scalebox{0.5}{
\begin{tikzpicture}[tensornetwork]
        \node[widetensor = 4] (A) at (0,0) {};

        \coordinate (A1) at (-1.5*\tensorsize, 0) {};
        \coordinate (A2) at (-0.5*\tensorsize, 0) {};
        \coordinate (A3) at (0.5*\tensorsize, 0) {};
        \coordinate (An) at (1.5*\tensorsize, 0) {};
        
        \draw (A.south -| A1) -- +(0, -0.5) node[below] {};
        \draw (A.south -| A2) -- +(0, -0.5) node[below] {};
        \draw (A.south -| A3)  +(0, -0.25) node[below] {$\cdots$};
        \draw (A.south -| An) -- +(0, -0.5) node[below] {};
       
\end{tikzpicture} 
}}=
    \vcenter{\hbox{\rotatebox{270}{\scalebox{0.35}{
\begin{tikzpicture}[tensornetwork]
    \filldraw[fill=Accent1!40, draw=none, rounded corners=10pt] (-0.7,1.5) rectangle (0.45,-6); 

    \filldraw[fill=Accent1!40, draw=none, rounded corners=10pt] (0.55 ,1.5) rectangle (1.7,-6);

    \filldraw[fill=Accent1!40, draw=none, rounded corners=10pt] (1.8,1.5) rectangle (2.95,-6); 
    
    \node[core] (A1) at (0,0) {};
    \node[core] (A2) at (0,-1.5) {};
    \node[nocore] (A3) at (0,-3) {$\vdots$};
    \node[core] (An) at (0,-4.5) {};

    \draw (A1) -- (A2) -- (A3) -- (An);
    \draw (A1.north west) -- + (-0.25,0.25) node[left] {};
    \draw (A2.north west) -- + (-0.25,0.25) node[left] {};
    \draw (An.north west) -- + (-0.25,0.25) node[left] {};

    \node[core] (B1) at (1.25,-4.5) {};
    \node[core] (B2) at (1.25,-3) {};
    \node[nocore] (B3) at (1.25,-1.5) {$\vdots$};
    \node[core] (Bn) at (1.25,0) {};
    
    \draw [rounded corners] (An.south) -- + (0,-0.5) -- + (1.25,-0.5) -- (B1.south);

    \draw (B1) -- (B2) -- (B3) -- (Bn);
    \draw (B1.north west) -- + (-0.25,0.25) node[left] {};
    \draw (B2.north west) -- + (-0.25,0.25) node[left] {};
    \draw (Bn.north west) -- + (-0.25,0.25) node[left] {};

    \node[core] (C1) at (2.5,0) {};
    \node[core] (C2) at (2.5,-1.5) {};
    \node[nocore] (C3) at (2.5,-3) {$\vdots$};
    \node[core] (Cn) at (2.5,-4.5) {};

    \node[nocore] (D) at (1.875,1.05) {$\cdots$};

    \draw [rounded corners] (Bn.north) -- + (0,0) |- (D) -- + (0.625,0) -- (C1.north);

    \draw (C1) -- (C2) -- (C3) -- (Cn);
    \draw (C1.north west) -- + (-0.25,0.25) node[right] {};
    \draw (C2.north west) -- + (-0.25,0.25) node[right] {};
    \draw (Cn.north west) -- + (-0.25,0.25) node[right] {};

\end{tikzpicture} 
}}}}\,.
    \label{eq: quantized tensor train}
\end{equation}

The QTT-format further reduces the storage cost to $\mathcal{O}(dlr^2) = \mathcal{O}(d\log(n) r^2)$, meaning that logarithmic complexity with respect to the mode sizes is achieved. Utilizing the same folding strategy on both indices of linear operators, allows them to be represented as QTT-operators (QTTO) \cite{kazeev_low-rank_2012}. This will prove to be useful for representing derivative operators in low-rank form.

\subsection{Update Equations}
The novel scheme is based on the conventional Yee FDTD technique \cite{yee_numerical_1966,taflove_computational_2005}. It employs a staggered discretization in both space and time and approximates the derivatives by second-order accurate central differences. In the spirit of finite integration (FIT), the scheme employs discrete grid voltages, i.e. fields rescaled by their respective edges \cite{weiland_time_1996}. Moreover, the electric field is rescaled with the free space impedance $Z_0$ and natural units are introduced for time, as such eliminating $\varepsilon_0$ and $\mu_0$ from the equations. The update equations for the $x$-components of the voltages are given by
\begin{equation}
    \Hn{n+{\scriptscriptstyle\tfrac{1}{2}}}{x} = \mathcal{R}\biggl[\Hn{n-{\scriptscriptstyle\tfrac{1}{2}}}{x} + \TENS{C}^h_x\odot\left(\TENSOP{D}_y (\En{n}{z})-\TENSOP{D}_z(\En{n}{y})\right)\biggr]\label{eq: faraday TT}
\end{equation}
\begin{equation}
    \En{n+1}{x} = \mathcal{R}\biggl[\En{n}{x}+  \TENS{C}^e_x\odot\left(\TENSOP{D}^\star_y ( \Hn{n+{\scriptscriptstyle\tfrac{1}{2}}}{z})-\TENSOP{D}^\star_z(\Hn{n+{\scriptscriptstyle\tfrac{1}{2}}}{y})\right)\biggr]
    ,\label{eq: ampere TT}
\end{equation}
where $\odot$ denotes the element-wise (or Hadamard) product. All objects in these equations need to be interpreted as being represented in the QTT-format, i.e., as QTTs or QTTOs. The binary operations ($+,-,\odot$) can be efficiently implemented in the QTT-format \cite{oseledets_tensor-train_2011}, see \cref{tab: comp complex op sat-qtt-fdtd}. The scaling is further discussed in \cref{sec: scaling}.

 The grid voltages are stored in the tensors $\Hn{n}{x}$ and $\En{n}{x}$. The operators $\TENSOP{D}$ and $\TENSOP{D}^\star$ represent the discretized derivative operators, on the primary and dual grid, respectively. These admit a QTTO decomposition with low rank, which depends on the boundary conditions employed; this is guaranteed by Lemma 3.1 in \cite{kazeev_low-rank_2012}.  In the case of a perfect electric conductor (PEC) at one side of a dimension and a perfect magnetic conductor (PMC) at the other end or periodic boundary conditions $r\leq2$ is obtained. In the case PECs or PMCs are enforced at both ends, the rank becomes $r\leq3$. Thus, the application of a derivative operator to a tensor $\TENS{A}$ represented in QTT-format results in a rank of $2r_A/3r_A$. The curl coefficients $\TENS{C}^e_x$ and $\TENS{C}^h_x$ are given by
\begin{equation}
        \TENS{C}^h_x = \Delta\tau\frac{\Delta x}{\Delta y\Delta z}[\bar{\boldsymbol{\mu}}^x_r]^{-1},
        \quad
        \TENS{C}^e_x = \Delta\tau\frac{\Delta x}{\Delta y\Delta z}[\bar{\boldsymbol{\varepsilon}}^x_r]^{-1},
\end{equation}
where, for simplicity, uniform gridding along the three spatial dimensions with grid steps $\Delta x$, $\Delta y$ and $\Delta z$ is assumed. $\Delta \tau$ denotes the time step in natural units.
The tensors $[\bar{\boldsymbol{\mu}}^x_r]$ and $[\bar{\boldsymbol{\varepsilon}}^x_r]$ contain the properly averaged material properties  and $^{-1}$ denotes the Hadamard inverse. As the Hadamard inverse of $[\bar{\boldsymbol{\mu}}^x_r]$ and $[\bar{\boldsymbol{\varepsilon}}^x_r]$ cannot be calculated explicitly in QTT-format, it is preferable to directly derive a QTT form for the inverse.
\begin{table}[!tb]
    \centering
    \caption{Computational complexity and resulting rank of the operations present in the FDTD scheme for FG and QTT.}
    \label{tab: comp complex op sat-qtt-fdtd}
    \begin{tabular}{cccc}
        \toprule
         Operation & FG & QTT&Resulting rank\\
         \midrule
         $\alpha\TENS{A}$& $\mathcal{O}(n^3)$ & $\mathcal{O}(\log(n)r_A)$ & $r_A$\\
         $\TENS{A}+\TENS{B}$&$\mathcal{O}(n^3)$&$\mathcal{O}(\log(n)(r_A+r_B)^2)$&$r_A+r_B$\\
         $\TENS{A}\odot\TENS{B}$&$\mathcal{O}(n^3)$&$\mathcal{O}(\log(n)r_A^2r_B^2)$&$r_Ar_B$\\
         $\TENSOP{D}|^{d'}_d(\TENS{A})$&$\mathcal{O}(n^3)$&$\mathcal{O}(\log(n)r_A^2)$&$2r_A$ / $3r_A$\\
         $\mathcal{R}(\TENS{A})$&/&$\mathcal{O}(\log(n)r_A^3)$&$r\leq r_A$\\
         \bottomrule
    \end{tabular}
\end{table}

Unfortunately, as becomes clear from \cref{tab: comp complex op sat-qtt-fdtd}, the rank increases when operations are performed. To mitigate the otherwise ever-increasing rank, a rounding routine that reduces the rank is required. To round the QTTs to a desired target rank $r$, the density matrix compression (DMC) algorithm (see \cite{mcculloch_density-matrix_2007,ye_quantized_2024}) is employed, which is denoted by $\mathcal{R}$.

The fixed-rank nature of this rounding procedure is different from the accuracy-based strategies of other low-rank solvers found in literature \cite{manzini_tensor-train_2023,scherzer_rank-limiting_2026,zhou_tensor-train_2025,nguyen_tensor_2026}. In situations dominated by numerical noise, such as when certain field components are identically zero in the exact solution, its accumulation can cause the ranks to grow over time. Fixed-rank rounding ensures ranks remain low, thereby enabling more predictable computational complexity. Alternative strategies have been proposed to mitigate unwanted rank growth such as regularization \cite{zhou_tensor-train_2025}, the group-rounding algorithm, accuracy-based rounding with rank caps \cite{scherzer_rank-limiting_2026} and smoothed contrast functions \cite{nguyen_tensor_2026}. 

Rounding is applied after a half-step has been completed, in order to preserve the full-rank structure of the update equations. Time stepping of this kind is called step-and-truncate (SAT) time integration \cite{einkemmer_review_2025,ye_practical_2026}. Hence, the advocated scheme is named the Step-and-Truncate Quantized Tensor Train Finite-Difference Time-Domain scheme or abbreviated SAT-QTT-FDTD.
\subsection{Current Sources}
To excite the computational domain, current sources are employed. To this end, the current is discretized on the same spatial grid as the electric voltages, while the same temporal discretization as the magnetic voltages is used. Subsequently, this results in an additional term in the update equations of the electric voltages, which is applied after the rounding of the electric voltage components has been performed, i.e., for the $x$-component
\begin{equation}
    \En{n+1}{x} = \left\{\En{n+1}{x}\right\} - \TENS{C}^j_x\odot \TENS{J}^{n+{\scriptscriptstyle\tfrac{1}{2}}}_x.
\end{equation}
In this equation, $\left\{\En{n+1}{x}\right\}$ denotes the right hand side of \cref{eq: ampere TT} and the current coefficient is given by
\begin{equation}
    \TENS{C}^j_x = \Delta\tau\Delta x[\bar{\boldsymbol{\varepsilon}}^x_r]^{-1}.
\end{equation}
In the case of a Hertzian dipole, utilized in the numerical experiment, $\TENS{C}^j_x \odot \TENS{J}^n_x$ admits a rank one QTT-decomposition.

\subsection{Perfectly Matched Layer}
To model open domains, a parabolically graded UPML is employed \cite{gedney_anisotropic_1996} and implemented using auxiliary differential equations (ADE). To this end, each grid voltage is complemented by one auxiliary degree of freedom, denoted by a dot, with corresponding additional update equation. For the $x$-component of the magnetic voltages, this results in
\begin{align}
    &\HDotn{n+{\scriptscriptstyle\tfrac{1}{2}}}{x} \!=\! \mathcal{R}\biggl[\TENS{S}^{\dot{h}}_x\!\odot\!\HDotn{n-{\scriptscriptstyle\tfrac{1}{2}}}{x}\!+\!\TENS{C}^{\dot{h}}_x\!\odot\!\left(\TENSOP{D}_y ( \En{n}{z})\!-\!\TENSOP{D}_z(\En{n}{y})\right)\biggr]\\
    &\Hn{n+{\scriptscriptstyle\tfrac{1}{2}}}{x} \!=\! \mathcal{R}\biggl[\TENS{S}^{h}_x \!\odot\! \Hn{n-{\scriptscriptstyle\tfrac{1}{2}}}{x} \!+ \!\TENS{A}^{h}_{x} \!\odot\!\HDotn{n+{\scriptscriptstyle \tfrac{1}{2}}}{x}\!+ \!\TENS{B}^{h}_{x} \!\odot\!\HDotn{n-{\scriptscriptstyle\tfrac{1}{2}}}{x}\biggr],\label{eq: PML}
\end{align}
with prefactors given by
\begin{equation}
    \TENS{S}^{\dot{h}}_x = \beta^-_y\odot(\beta_y^+)^{-1}\;,\;\TENS{S}^h_x = \beta_z^-\odot(\beta_z^+)^{-1}\;,
\end{equation}
\begin{equation}
    \TENS{C}^{\dot{h}}_x = (\beta_y^+)^{-1} \odot \TENS{C}^{h}_x\;,
\end{equation}
\begin{equation}
    \TENS{A}^{h}_{x} = \beta_x^+ \odot (\beta_z^+)^{-1}\;,\; \TENS{B}^{h}_{x} = -\beta_x^- \odot (\beta_z^+)^{-1}.
\end{equation}
The tensors $\beta_d^\pm$ are given by $
    \beta_d^{\pm} = \frac{\kappa_d}{\Delta\tau} \pm \frac{Z_0\sigma_d}{2}\, , \, d\in\{x,y,z\}$,
were $\kappa_d$ and $\sigma_d$ are the tensors containing the real and imaginary stretching parameters of the PML, respectively.
The electric and remaining magnetic update equations are analogous. The prefactors are represented in QTT-format, via TT-SVD, with an accuracy $\varepsilon =$ 10$^{-10}$, resulting in a representation with moderate ranks, more specifically $r_k \leq 8$.
\subsection{Scaling}\label{sec: scaling}
The scaling of the numerical operations present in the update equations is given in \cref{tab: comp complex op sat-qtt-fdtd} and compared to the FG case. From the table it becomes clear that the rounding is the most expensive part of the framework. One could argue that the limiting step would be the Hadamard product, because of the $r_A^2r_B^2$ factor, which appears to scale with power four. However, it scales only quadratically with respect to the rank of the fields; the rank of the curl coefficients is fixed for a certain simulation.
Moreover, the rank of the QTTs that are rounded is larger than the ranks of the QTTs in the binary operations, as it is the last step in an update. 
As a result, the cost per time step of the scheme scales as $\mathcal{O}(\log(n)r^3)$. Provided that the ranks remain moderate, the SAT-QTT-FDTD scheme improves the volumetric scaling into logarithmic scaling. The memory cost exhibits $\mathcal{O}(\log(n)r^2)$ scaling, as explained in \cref{sec: QTT}.
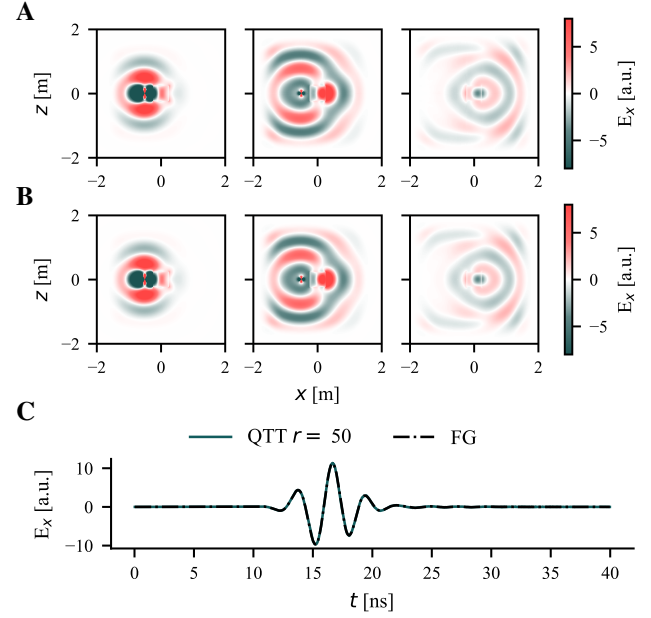
\begin{figure}[!tb]
    \centering
    \begin{tikzpicture}
        \node[anchor=north west] at (0,0) {\includegraphics[width=0.95\columnwidth]{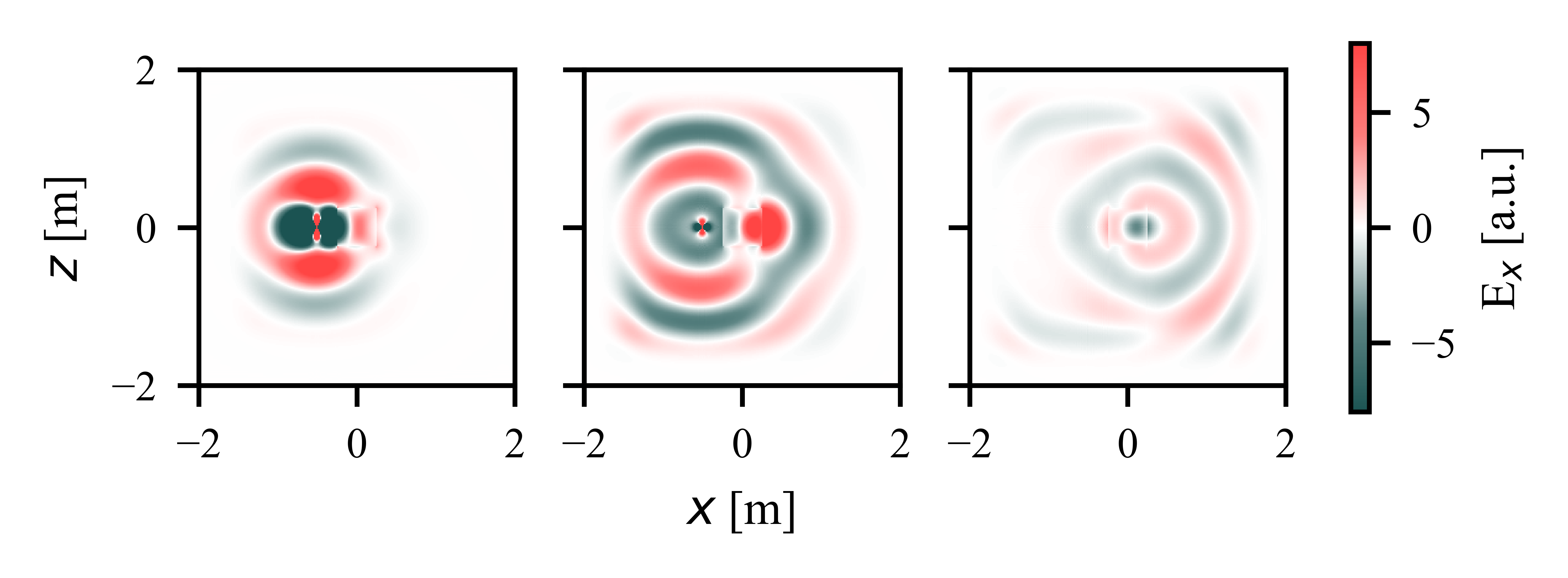}};
        \node[anchor=north west] at (0,0) {\textbf{A}};
    \end{tikzpicture}
    \vspace{-2.5em}
    
    \begin{tikzpicture}
        \node[anchor=north west] at (0,0) {\includegraphics[width=0.95\columnwidth]{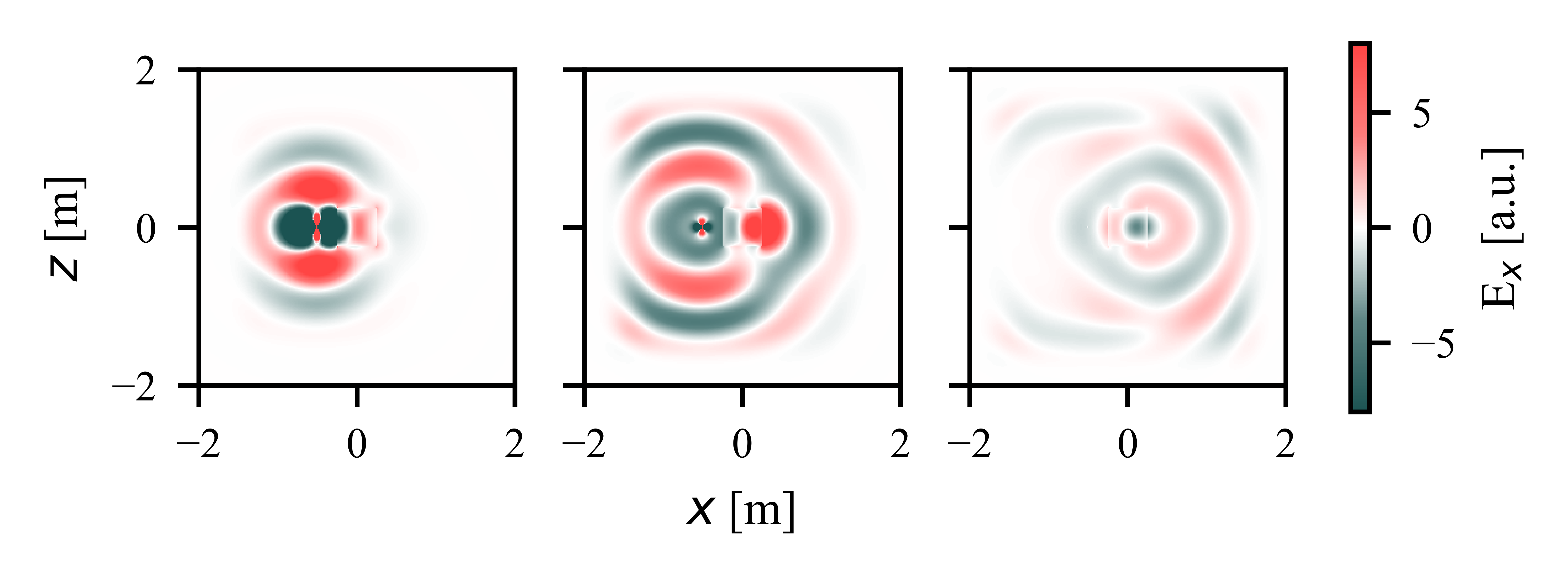}};
        \node[anchor=north west] at (0,0) {\textbf{B}};
    \end{tikzpicture}
    \vspace{-1.4em}
    
    \begin{tikzpicture}
        \node[anchor=north west] at (0,0) {\includegraphics[width=0.95\columnwidth]{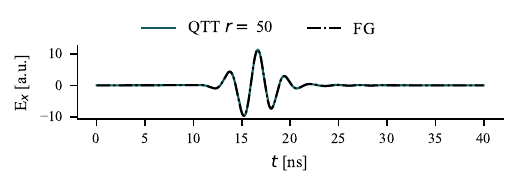}};
        \node[anchor=north west] at (0,0) {\textbf{C}};
    \end{tikzpicture}
    \vspace{-1.5em}
    \caption{Snapshots of the $x$-component of the electric field at $y=0$, for three instances: \SI{12.50}{\nano\second}, \SI{16.25}{\nano\second} and \SI{20}{\nano\second}. (A) SAT-QTT-FDTD $r=50$. (B) FG FDTD. (C) The $x$-component of the electric field at (\SI{0.5}{\meter}, \SI{0}{\meter}, \SI{0}{\meter}) as a function of time.}
    \label{fig: snapshots}
\end{figure}

\section{Numerical Results}\label{sec:numerical results}
As a validation example a dielectric cube with $\varepsilon_r =$ 4 and edges of \SI{50}{\centi\meter} is illuminated by a Hertzian dipole oriented normal to one of the cube's sides. The cube's center is placed at the origin and the dipole is located at (-\SI{0.5}{\meter}, \SI{0}{\meter}, \SI{0}{\meter}). To mimic open boundaries, the simulation domain $\Omega = [-\SI{1}{\meter} , \SI{1}{\meter}]^3$ is surrounded by a \SI{1}{\meter} thick PML.  A Gaussian modulated sine is used as the source function of the excitation, i.e., 
\begin{equation}
    s\left(t\right) = Ae^{{-\frac{\pi^2\Delta f_\text{FWHM}^2(t-t_0)^2}{4\ln2}}}\sin\left(2\pi f_c(t-t_0)\right),
\end{equation}
where $f_c = \SI{300}{\mega\hertz}$, $\Delta f_\text{FWHM}=\SI{200}{\mega\hertz}$ and $t_0 = \SI{12}{\nano\second}$.
The time step is chosen such that the Courant number $\text{CN} = 0.9$.
In the case of a cube, $[\bar{\boldsymbol{\varepsilon}}^d_r]^{-1}$ admits a QTT-decomposition with rank theoretically bound from above by 5. This follows from the fact that Heaviside functions have QTT rank $r=2$ \cite{gourianov_exploiting_2022}. In practice lower ranks $r_k\leq$ 3 were obtained using TT-SVD with an accuracy set to 10$^{-8}$.
The experiments were performed on a machine equipped with two Intel Xeon Gold 6226, 504 Gb RAM and running Ubuntu 24.04 using an in-house developed python code.

Snapshots of the $x$-component of the electric field at \SI{12.5}{\nano\second}, \SI{16.25}{\nano\second} and \SI{20}{\nano\second} calculated on a $256\times 256\times 256$ grid are shown in \cref{fig: snapshots} (A) and (B) for SAT-QTT-FDTD at rank 50 and traditional FG FDTD, respectively. The scattering at the cube is clearly visible in both the QTT and FG snapshots. Moreover, the QTT snapshots show excellent correspondence with the FG snapshots. \cref{fig: snapshots}~(C) shows the time-evaluation of the $x$-component of the electric field at (\SI{0.5}{\meter}, \SI{0}{\meter}, \SI{0}{\meter}); the QTT results match the FG ones. 
Additionally the PML clearly absorbs the fields, highlighting a successful implementation. Below, quantitative results concerning the accuracy-efficiency trade-off are presented.

\begin{figure}[!tb]
    \centering
    \vspace{-0.6em}
    \includegraphics{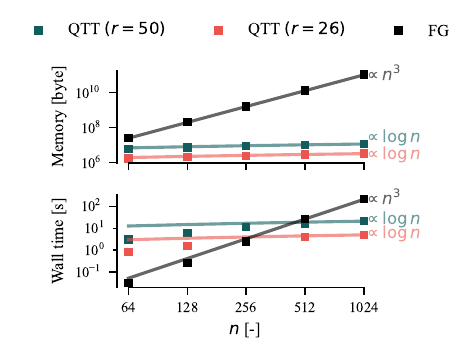}
    \vspace{-1.6em}
    \caption{Experimentally observed scaling (squares) for FG FDTD and SAT-QTT-FDTD compared to the theoretically expected scaling (lines).}
    \label{fig:scaling}
\end{figure}

To evaluate the complexity of the method, with respect to the number of grid points along a certain dimension $n=2^l$, 100 time steps were performed on different $n\times n\times n$ grids and the average cost per time step was calculated. For SAT-QTT-FDTD the experiment was conducted at two different ranks ($r=26,50$). The results are depicted in \cref{fig:scaling}. The observed scalings of the memory requirements and of the wall time per time step correspond to the theoretically expected ones, namely $\mathcal{O}(n^3)$ and $\mathcal{O}(\log n)$ for FG FDTD and SAT-QTT-FDTD, respectively. In the experiment, SAT-QTT-FDTD attains significant memory compression ratios at every discretisation. When $l = 9$, memory compression ratios of over three orders of magnitude are observed. When $l = 10$, a memory compression ratio exceeding four orders of magnitude is achieved for $r=26$. 

\begin{figure}[!tb]
    \centering
    \vspace{-0.6em}
    \includegraphics{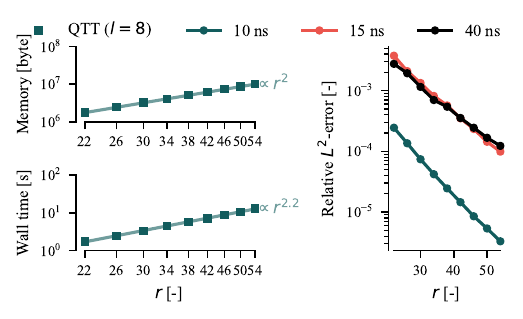}
    \vspace{-2.5em}
    \caption{Resources (left) and accuracy with respect to the FG solver (right) of the SAT-QTT-FDTD ($l=8$) method as a function of rank. Accuracy is evaluated after \SI{10}{\nano\second}, \SI{15}{\nano\second} and \SI{40}{\nano\second}.}
    \label{fig:scaling rank}
\end{figure}

In addition to the discretization, the rank constitutes a second parameter that controls both accuracy and required resources. To investigate its influence, the simulation is repeated for different ranks at a fixed discretization $l=8$. The accuracy is evaluated using the relative $L^2$-error with respect to a FG FDTD solver on the same grid, defined by
\begin{equation}
    \mathcal{E}^n_{L^2} = \frac{\sqrt{\sum_d\norm{\En{n,QTT}{d}-\En{n,FG}{d}}^2_{L^2}}}{\max_n\sqrt{\sum_d\norm{\En{n,FG}{d}}^2_{L^2}}} \, , \, d\in\{x,y,z\},
\end{equation}
in which $\norm{\cdot}_{L^2}$ represents the discrete $L^2$-norm. The wall time per time step is again evaluated by taking the average of 100 time steps. The results of this experiment are depicted in \cref{fig:scaling rank}. 
The memory exhibits the expected $\mathcal{O}(r^2)$ complexity. The observed scaling of the wall time per time step $\mathcal{O}(r^{2.2})$ is better than with the theoretical predicted one $\mathcal{O}(r^3)$. The slight deviation is attributed to the fact that $\mathcal{O}(d\log(n)r^3)$ is a worst case estimate for the DMC algorithm and the true scaling relation is more complex \cite{ye_quantized_2024}. The right panel of the figure shows the relative $L^2$-error at \SI{10}{\nano\second}, \SI{15}{\nano\second} and \SI{40}{\nano\second}, i.e., before $t_0$, after $t_0$ and at the end of the simulation, respectively. As the rank increases the error decreases, as expected. The error grows during the initial energy injection by the source. However, once the energy has been introduced into the domain, the error stabilizes. Excellent accuracy, with errors dropping considerably below $10^{-3}$, is observed.

\section{Conclusion}\label{sec:conclusion}
This letter proposes a QTT-accelerated FDTD framework for the simulation of 3-D open-space scattering problems. The validation example shows that the method achieves excellent accuracy with respect to FG FDTD, while reducing the cubic scaling associated with FDTD to a logarithmic scaling. The rank allows for a controllable trade-off between accuracy and computational resources. Owing to its superior scaling, the proposed method achieves memory savings of several orders of magnitude, illustrating its potential for the simulation of large-scale electromagnetic systems.
\bibliographystyle{IEEEtran}
\bibliography{LetterBibliography}
\end{document}